\documentclass[12pt]{article}

\usepackage{color}
\usepackage{latexsym,dsfont}
\usepackage{amssymb}
\usepackage{amsmath}
\usepackage{amsthm}
\usepackage{graphicx}
\usepackage{enumerate}
\usepackage{hyperref}

\newtheorem{Theorem}{Theorem}[part]
\newtheorem{Definition}{Definition}[part]
\newtheorem{Proposition}{Proposition}[part]

\newtheorem{Lemma}{Lemma}[part]

\newtheorem{Remark}{Remark}[part]

\def \R{\mathbb{R}}

\def \E{\mathbb{E}}
\def \F{\mathbb{F}}
\def \G{\mathbb{G}}

\def \P{\mathbb{P}}
\def \Q{\mathbb{Q}}
\def \D{\mathbb{D}}

\def \PMc{{\cal P \cal M}}

\def \Bc{{\cal B}}

\def \Dc{{\cal D}}
\def \Ec{{\cal E}}
\def \Fc{{\cal F}}
\def \Gc{{\cal G}}

\def \Pc{{\cal P}}

\def \Sc{{\cal S}}
\def \Tc{{\cal T}}

\def \1{\mathds{1}}
\def \XunPi{X^{1, \pi}}

\def \YunPi{Y^{1, \pi}}
\def \ZunPi{Z^{1, \pi}}

\def \YzePi{Y^{0, \pi}}
\def \ZzePi{Z^{0, \pi}}
\def \XzePi{X^{0, \pi}}

\def \ni{\noindent}

\def \eps{\varepsilon}

\def \ep{\hbox{ }\hfill$\Box$}

\def\reff#1{{\rm(\ref{#1})}}

\def\beqs{\begin{eqnarray*}}
\def\enqs{\end{eqnarray*}}
\def\beq{\begin{eqnarray}}
\def\enq{\end{eqnarray}}

\newcommand{\nc}{\newcommand}
\nc{\esssup}{\mathop{\mathrm{ess\;sup}}}

\begin{document}
\title{A decomposition approach for\\ the discrete-time approximation of \\BSDEs with a jump II: the quadratic case}
\author{Idris Kharroubi \footnote{The research of the author benefited from the support of the French ANR research grant LIQUIRISK} \\
\small{CEREMADE, CNRS  UMR 7534}\\
\small{Universit\'e Paris Dauphine} \\
\small{\texttt{kharroubi @ ceremade.dauphine.fr}}\and 
Thomas Lim \footnote{The research of the author benefited from the support of the ``Chaire Risque de Cr\'edit'', F\'ed\'eration Bancaire Fran\c caise}\\
\small{Laboratoire d'Analyse et Probabilit\'es,} \\
\small{Universit\'e d'Evry and ENSIIE,}\\
\small{\texttt{thomas.lim @ ensiie.fr} }}

\date{}

\maketitle
\begin{abstract}
We study the discrete-time approximation for solutions of quadratic forward backward stochastic differential equations (FBSDEs) driven by a Brownian motion and  a jump process which could be dependent. Assuming that the generator has a quadratic growth w.r.t. the variable z and  the terminal condition is bounded, we prove the convergence of the scheme when the number of time steps n goes to infinity. Our approach is based on the companion paper \cite{kl11} and allows to get a convergence rate similar to that of schemes of Brownian FBSDEs. 
\end{abstract}


\vspace{1cm}

\ni \textbf{Keywords:} discrete-time approximation, forward-backward SDE with a jump, generator of quadratic growth, progressive enlargement of filtrations, decomposition in the reference filtration. 

\vspace{1cm}

\ni\textbf{MSC classification (2000):} 65C99, 60J75, 60G57.

\newpage

\tableofcontents

\newpage

\setcounter{section}{0}

\section{Introduction}
\setcounter{equation}{0} \setcounter{Assumption}{0}
\setcounter{Theorem}{0} \setcounter{Proposition}{0}
\setcounter{Corollary}{0} \setcounter{Lemma}{0}
\setcounter{Definition}{0} \setcounter{Remark}{0}

In this paper, we study a discrete-time approximation for the solution of a forward-backward stochastic differential equation (FBSDE) with a jump and taking the following form
 \begin{equation*} \left\{
\begin{aligned}
  X_t ~ & = ~ x + \int_0^t b(s, X_s) ds + \int_0^t \sigma(s) dW_s + \int_0^t \beta(s, X_{s^-}) dH_s \;,\\
  Y_t ~ & = ~ g(X_T) + \int_t^T f(s, X_s, Y_s, Z_s, U_s) ds - \int_t^T Z_s dW_s - \int_t^T U_s dH_s\;, \\
  \end{aligned}
\right. 0\leq t\leq T,\end{equation*}
where $H_t = \1_{\tau \leq t}$ and $\tau$ is a jump time, depending on $W$, which can represent a default time in credit risk or counterparty risk. Such equations naturally appear in finance, see for example Bielecki and Jeanblanc \cite{biejea08}, Lim and Quenez \cite{limque09}, Ankirchner \emph{et al.} \cite{ankblaeyr09} for an application to the exponential utility maximization problem and Kharroubi and Lim \cite{kl11} for the hedging problem in a complete market. 
The approximation of such equation is therefore of important interest for practical applications in finance.

This paper is the second part of a two sides work studying the discretization of such FBSDEs. The first part \cite{kl11a} deals with the case of a Lipschitz continuous generator $f$ and this second part focuses on the case of  a generator $f$ with quadratic growth w.r.t. $z$.



In the literature, the problem of discretization of FBSDEs  with Lipschitz generator has been widely studied in the Brownian framework, \textit{i.e.} no jump, see e.g. \cite{maproyon94, doumapro96, che97, boutou04, 
zha04, delmen06}. 
More recently, the case of quadratic generators w.r.t. $z$ has been considered by Imkeller \emph{et al.} \cite{imkreizha10} and Richou \cite{ric10}.
 
The discretization of BSDEs with jumps has been studied by Bouchard and Elie \cite{boueli08} in the case where the  generator is  Lipschitz continuous and the jumps are independent of the Brownian motion. 

In this paper, we study the discrete-time approximation of FBSDEs with
\begin{itemize}
\item a generator satisfying a quadratic growth w.r.t. the variable  $z$,

\item  a jump time $\tau$ that could be dependent on the Brownian motion via the density assumption.

\end{itemize}

Such FBSDEs are of interest in finance since utility maximization problems lead in general to quadratic generators and default time models generally impose the dependence of $\tau$ w.r.t. the Brownian motion. 


To study the discretization of such equations, we can not directly work on the BSDE as done in  \cite{boueli08} since their approach uses 
 a regularity result for the process $Z$ based on Malliavin calculus. In our situation, the dependence of $\tau$ w.r.t. $W$ prevent us from using such an approach since no Malliavin calculus theory has been set for this framework. 
 
 To get a discrete-time approximation scheme we then use the results of \cite{kl11}, that allow to decompose the FBSDE with a jump in a recursive system of Brownian FBSDEs. 
 We then provide estimates on the solutions to each Brownian FBSDE.
These estimates allow to prove that the solution satisfies a  Lipschitz FBSDE. 
 Finally, we obtain a discrete-time approximation scheme by using the results of the previous part \cite{kl11a} on Lipschitz BSDEs. 

The paper is organized as follows. The next section presents the FBSDE, the different assumptions on the coefficients of the FBSDE and recalls the result of \cite{kl11}. In Section 3, we give some estimations on the solution of the FBSDE. In Section 4, we give a discrete-time approximation scheme for the FBSDE and provide a global error estimate.

\section{Preliminaries}
\setcounter{equation}{0} \setcounter{Assumption}{0}
\setcounter{Theorem}{0} \setcounter{Proposition}{0}
\setcounter{Corollary}{0} \setcounter{Lemma}{0}
\setcounter{Definition}{0} \setcounter{Remark}{0}

\subsection{Notation}
Throughout this paper, we let $(\Omega, \Gc, \P)$ a complete probability space on which is defined a standard one dimensional Brownian motion $W$. We denote 
 $\F :={(\Fc_t)}_{t\geq 0}$ the natural filtration of $W$ augmented by all the $\P$-null sets. We also consider on this space a random time $\tau$, i.e. a nonnegative $\Fc$-measurable random variable, and we denote classically the associated jump process by $H$ which is given by 
\beqs
H_t & := & \1_{\tau \leq t}\;,\quad t\geq 0\;. 
\enqs
We denote by  $\D := (\Dc_t)_{t \geq 0}$ the smallest right-continuous filtration for which $\tau$ is a stopping time. The global information is then defined by the progressive enlargement $\G := {(\Gc_t)}_{t \geq 0}$ of the initial filtration where \beqs
\Gc_t &  :=  & \bigcap_{\eps>0} \Big(\Fc_{t+\eps} \vee \Dc_{t+\eps}\Big)
\enqs
for all $t\geq 0$. 
 This kind of enlargement  was introduced by Jacod, Jeulin and Yor in the 80s (see e.g. \cite{jeu80}, \cite{jeuyor85} and \cite{jac87}). We introduce some notations used throughout the paper\begin{itemize}
\item $\Pc(\F)$ (resp. $\Pc(\G)$) is the $\sigma$-algebra of $\F$ (resp. $\G$)-predictable measurable subsets of $\Omega \times \R_{+}$, i.e. the $\sigma$-algebra generated by the left-continuous $\F$ (resp. $\G$)-adapted processes, 
\item $\PMc(\F)$ (resp. $\PMc(\G)$) is the $\sigma$-algebra of $\F$ (resp. $\G$)-progressively measurable subsets of $\Omega\times\R_{+}$. 
\end{itemize}
We shall make, throughout the sequel, the standing assumption in the progressive enlargement of filtrations known as  density assumption (see e.g. \cite{jiapha09, jkp11,kl11}).

\vspace{2mm}

\ni \textbf{(DH)} There exists a positive and bounded $\Pc(\F)\otimes\Bc(\R_+)$-measurable process $\gamma$ such that
\beqs
\P \big[\tau\in d\theta \;\big| \;\Fc_t \big] & = & \gamma_t(\theta)d\theta \;,\quad t\geq 0\;.
\enqs

\vspace{2mm}

Using Proposition 2.1 in \cite{kl11} we get that \textbf{(DH)}  ensures that the process $H$ admits an intensity.
 \begin{Proposition}
The process $H$ admits a compensator of the form $\lambda_tdt$, where the process $\lambda$ is defined by 
\beqs
\lambda_t & := & \frac{\gamma_t(t)}{\P \big[\tau>t\; \big| \;\Fc_t\big]}\mathds{1}_{t\leq \tau}\;,\quad  t\geq 0\;.
\enqs
\end{Proposition} 
\ni We impose the following assumption to the process $\lambda$.

\vspace{2mm}

\ni\textbf{(HBI)} The process $\lambda$ is bounded. 

\vspace{2mm}

\ni We also introduce the martingale invariance assumption known as the \textbf{(H)}-hypothesis.

\vspace{2mm}

\ni \textbf{(H)} Any $\F$-martingale remains a $\G$-martingale.

\vspace{2mm}

 
\ni We now introduce the following spaces, where $a,b\in\R_+$ with $a< b$,
and $T < \infty$ is the terminal time.
\begin{itemize}
\item $\Sc_{\G}^\infty[a,b]$ (resp. $\Sc_{\F}^\infty[a,b]$) is the set of $\PMc(\G)$ (resp.  $\PMc(\F)$)-measurable processes $(Y_t)_{t\in[a,b]}$ essentially bounded
\beqs
{\| Y \|}_{\Sc^\infty[a,b]} & := & \esssup_{t\in[a,b]}|Y_{t}|~<~\infty \;.
\enqs
\item $\Sc_{\G}^p[a,b]$ (resp. $\Sc_{\F}^p[a,b]$), with $p\geq 2$, is the set of $\PMc(\G)$ (resp.  $\PMc(\F)$)-measurable processes $(Y_t)_{t\in[a,b]}$ such that
\beqs
{\| Y \|}_{\Sc^p[a,b]} & := & \Big(\E\Big[\sup_{t\in[a,b]}|Y_{t}|^p\Big]\Big)^{1\over p}~<~\infty \;.
\enqs

\item $H^p_{\G}[a,b]$ (resp. $H^p_{\F}[a,b]$), with $p\geq 2$, is the set of $\Pc(\G)$ (resp. $\Pc(\F)$)-measurable processes $(Z_t)_{t\in[a,b]}$ such that
\beqs
{\|Z\|}_{H^p[a,b]}  & := & \E\Big[ \Big(\int_a^b |Z_t|^2 dt \Big)^{p\over 2}\Big]^{1\over p} ~< ~\infty \;.
\enqs
\item $L^2(\lambda)$ is the set of $\Pc(\G)$-measurable processes $(U_t)_{t \in [0, T]}$ such that
\beqs
{\|U\|}_{L^2(\mu)} &:= & \Big(\E\Big[\int_0^T \lambda_s |U_{s}|^2ds\Big]\Big)^{1\over2}~<~\infty\;.
\enqs

\end{itemize}

  \subsection{ Quadratic growth Forward-Backward SDE with a jump}
Given measurable functions $b:[0,T]\times\R \rightarrow \R$, $\sigma:[0,T] \rightarrow \R$, $\beta: [0,T] \times \R \rightarrow\R$, $g: \R \rightarrow\R$ and $f: [0, T] \times \R \times\R\times\R \times \R\rightarrow\R$, and an initial condition $x \in \R$, we  study the discrete-time approximation of  the solution $(X, Y, Z, U)$ in $\Sc^2_\G [0, T] \times \Sc^\infty_\G [0, T] \times H^2_\G [0, T] \times L^2(\lambda)$ to the following forward-backward stochastic differential equation
\begin{eqnarray}\label{FSDE} 
  X_t  & = & x + \int_0^t b(s, X_s) ds + \int_0^t \sigma(s) dW_s + \int_0^t \beta(s, X_{s^-} ) dH_s \;, \quad 0 \leq t\leq T\;,\\ \nonumber
  Y_t  & =  & g(X_T) + \int_t^T f\big(s, X_s, Y_s, Z_s, U_s(1-H_s)\big) ds\\ \label{BSDE}
   & & \quad \quad \quad \quad  - \int_t^T Z_s dW_s - \int_t^T U_s dH_s\;, \quad 0 \leq t\leq T\;, ~~
\end{eqnarray}
  when the generator of the BSDE \reff{BSDE} has a quadratic growth w.r.t. $z$. 
  \begin{Remark}
  {\rm In the BSDE \reff{BSDE}, the jump component $U$ of the unknown $(Y,Z,U)$ appears in the generator $f$ with an additional multiplicative term $1-H$. This ensures the equation to be well posed in $\Sc_\G^\infty[0,T]\times H^2_\G[0,T]\times L^2(\lambda)$. Indeed,  the component $U$ lives in $L^2(\lambda)$, thus its value on $(\tau\wedge T,T]$ is not defined, since the intensity $\lambda$ vanishes on $(\tau\wedge T,T]$. We therefore introduce the term $1-H$ to kill the value of $U$ on $(\tau\wedge T,T]$ and hence to avoid making the equation depends on it.   }
  \end{Remark}
  We first prove that the decoupled  system \reff{FSDE}-\reff{BSDE} admits a solution.
  To this end, we introduce several  assumptions on the coefficients  $b$, $\sigma$, $\beta$, $g$ and $f$. 
  We consider the following assumptions for the forward coefficients.
  
\vspace{2mm}  
  
 \ni\textbf{(HF)} There exist two constants $K_a$ and $L_a$ such that 
  the functions $b$, $\sigma$ and $\beta$ satisfy 
  \beqs
  |b(t, 0)|+ |\sigma (t)|+|\beta(t,0)| & \leq  & K_a \;,
  \enqs
and
  \beqs
  |b(t, x) - b(t, x')|  +  |\beta (t, x) - \beta(t, x')|  & \leq & L_a |x - x'| \;,
  \enqs
  for all $(t,x,x')\in[0,T]\times\R \times \R$. 
  
\vspace{2mm}  
  

  

\ni For the backward coefficients $g$ and $f$, we consider the following assumptions.
  
\vspace{2mm}  
  
%
%
%
%

  

 \ni\textbf{(HBQ)} \begin{itemize}
 \item There exist three constants $M_g$, $K_g$ and $K_q$ such that 
  the functions $g$ and $f$ satisfy 
  \beqs
  |g(x)| & \leq & M_g\;,\\
  |g(x)-g(x')| & \leq & K_g |x-x'| \;,\\
  |f(t,x,y,z,u) - f(t,y',z,u)| & \leq & K_q |y-y'| \;,\\
  |f (t,x,y,z,u)|  & \leq & 
  K_q \big(1+ |y|+|z|^2+|u|\big) \;,
  \enqs
 for all $(t,x,x',y,y',z,u)\in [0,T] \times \R^2 \times \R^2 \times \R \times\R$. 
 \item For any $R >0$ there exists a function $mc^f_R$ such that $\lim_{\epsilon \rightarrow 0} mc^f_R(\epsilon) =0$ and
 \beqs
 |f(t,x,y,z, u-y) - f(t,x,y',z', u-y') | & \leq & mc^f_R(\epsilon)
 \enqs
 for all $(t,x,y,y',z,z',u) \in [0,T] \times\R \times  \R^2 \times \R^2 \times \R$ s.t. $|y|$, $|z|$, $|y'|$, $|z'| \leq R$ and $|y-y'| + |z-z'| \leq \epsilon$. 
 \item $f(t,.,u)=f(t,.,0)$ for all $u \in \R$ and all $t \in (\tau \wedge T, T]$.
 \item The function $f(t,x,y,.,u)$ is convexe (or concave) uniformly in  $(t,x,y,u)\in [0,T] \times \R
  \times\R\times\R$.  \end{itemize}
%

\vspace{2mm}  
\ni In the sequel $K$ denotes a generic constant appearing in \textbf{(HBQ)} and \textbf{(HF)} and which may vary from line to line.\\


\ni In the purpose to prove the existence of a solution to the FBSDE (\ref{FSDE})-(\ref{BSDE}) we follow the decomposition approach initiated by \cite{kl11} and for that we introduce the recursive system of  FBSDEs associated with 
\reff{FSDE}-\reff{BSDE}.

\vspace{2mm}

\ni$\bullet$ Find $(X^1(\theta),Y^1(\theta),Z^1(\theta))\in\Sc^2_\F[0,T]\times\Sc_\F^\infty[\theta ,T]\times H_\F^2[\theta,T]$ such that 
\beq\label{FSDE1} 
  X^1_t(\theta)  & = & x + \int_0^t b\big(s, X^1_s(\theta)\big) ds + \int_0^t \sigma(s) dW_s + \beta\big(\theta, X^1_{\theta^-} (\theta)\big) \1_{\theta \leq t} \;, \quad 0\leq t\leq T\;,\\
\label{BSDE1}  Y^1_t\big(\theta\big)  & = & g\big(X^1_T(\theta)\big) + \int_t^T f\big( s, X^1_s(\theta), Y^1_s(\theta), Z^1_s(\theta), 0\big) ds - \int_t^T Z^1_s(\theta) dW_s \;, \quad \theta\leq t\leq T\;,~\quad\quad 
\enq
for all $\theta\in[0,T]$.

\vspace{3mm}

\ni$\bullet$ Find $(X^0,Y^0,Z^0)\in\Sc^2_\F[0,T]\times\Sc_\F^\infty[0 ,T]\times H_\F^2[0,T]$ such that 
\beq\label{FSDE0}
  X^0_t  & = & x + \int_0^t b\big(s, X^0_s\big) ds + \int_0^t \sigma(s ) dW_s \;, \quad 0\leq t\leq T\;,\\ \label{BSDE0}
  Y^0_t  & = & g\big(X^0_T\big) + \int_t^T f \big( s, X^0_s, Y^0_s, Z^0_s, Y^1_s(s) - Y^0_s\big) ds - \int_t^T Z^0_s dW_s \;, \quad 0\leq t\leq T\;.\quad
\enq
Then, the link between the FBSDE \reff{FSDE}-\reff{BSDE} and the recursive system of FBSDEs \reff{FSDE1}-\reff{BSDE1} and \reff{FSDE0}-\reff{BSDE0} is given by the following result.
\begin{Theorem}  \label{theoreme existence unicite}
Assume that \textbf{(DH)}, \textbf{(HBI)}, \textbf{(H)}, \textbf{(HF)} and \textbf{(HBQ)} 
 hold true. Then, the FBSDE \reff{FSDE}-\reff{BSDE} admits a unique solution $(X,Y,Z,U)\in\Sc^2_\G[0,T]\times\Sc^\infty_\G[0,T]\times H^2_\G[0,T]\times L^2(\lambda)$ given by
\begin{equation}\label{expression solution} \left\{
\begin{aligned}
X_t ~ & = ~ X^0_t \1_{t < \tau} + X^1_t(\tau) \1_{\tau \leq  t} \;,\\
Y_{t}  ~& =~ Y^0_t \mathds{1}_{t< \tau} + Y^1_t(\tau) \mathds{1}_{\tau \leq t} \;,\\
Z_{t} ~& =  ~Z^0_t \mathds{1}_{t \leq \tau} + Z^1_t(\tau) \mathds{1}_{\tau <  t} \;,\\
U_{t}~ & = ~ \big(Y^1_t(t) -Y^0_t\big) \mathds{1}_{t \leq \tau} \;,
\end{aligned}
\right.\end{equation}
 where  $(X^1(\theta),Y^1(\theta),Z^1(\theta))$ is the unique solution to the FBSDE \reff{FSDE1}-\reff{BSDE1} in $\Sc^2_\F[0,T]\times\Sc^\infty_\F[\theta,T]\times H^2_\F[\theta,T]$, for $\theta\in[0,T]$, and $(X^0,Y^0,Z^0)$ is the unique solution to the FBSDE \reff{FSDE0}-\reff{BSDE0} in $\Sc^2_\F[0,T]\times\Sc^\infty_\F[0,T]\times H^2_\F[0,T]$. 
\end{Theorem}
\ni\textbf{Proof.}
%
%
The existence and uniqueness of the forward process $X$ and its link with $X^0$ and $X^1$ have already been proved in the first part of this work \cite{kl11a}. We now concentrate on the backward equation. 


\ni  To follow the decomposition approach initiated by the authors in \cite{kl11}, we need the generator to be predictable. To this end, we notice that in the BSDE \reff{BSDE}, we can replace the generator $(t,y,z,u)\mapsto f(t,X_t,y,z,u(1-H_t))$ by the predictable map $(t,y,z,u)\mapsto f(t,X_{t^-},y,z,u(1-H_{t^-}))$. 

\ni Suppose that \textbf{(DH)}, \textbf{(HBI)}, \textbf{(H)} and \textbf{(HBQ)} hold true. 
The existence of a solution $(Y,Z,U)\in\Sc^2_\G[0,T]\times H^2_\G[0,T]\times L^2(\lambda)$ is  then a direct consequence of Proposition 3.1 in \cite{kl11}.
We then notice that from the definition of $H$ we have $f(t,x,y,z,u(1-H_t))=f(t,x,y,z,0)$ for all $t\in(\tau\wedge T,T]$.  This property and \textbf{(DH)}, \textbf{(HBI)}, \textbf{(H)} and \textbf{(HBQ)} 
 allow to apply  Theorem 4.2 in \cite{kl11}, which gives the uniqueness of a solution to BSDE \reff{BSDE}.   
\ep

\vspace{2mm}

Throughout the sequel, we give an approximation of the solution to the FBSDE \reff{FSDE}-\reff{BSDE} by studying the approximation of the solutions to the recursive system of FBSDEs \reff{FSDE1}-\reff{BSDE1} and \reff{FSDE0}-\reff{BSDE0}.

\section{A priori estimation on the gain process} 

\setcounter{equation}{0} \setcounter{Assumption}{0}
\setcounter{Theorem}{0} \setcounter{Proposition}{0}
\setcounter{Corollary}{0} \setcounter{Lemma}{0}
\setcounter{Definition}{0} \setcounter{Remark}{0}

Before giving the discrete-time scheme for the FBSDE (\ref{FSDE})-(\ref{BSDE}) we give a uniform bound for the processes $Z^0$ and $Z^1$ which allows to prove that the BSDE \reff{BSDE} is Lipschitz and thus we can use the discrete-time scheme given in \cite{kl11a}. For that we introduce the BMO-martingales class, and we also give some bounds for the processes $X^0$, $X^1$, $Y^0$ and $Y^1$.

\subsection{BMO property for the solution of the BSDE}

To obtain a uniform bound for the processes $Z^0$ and $Z^1$ we need the following assumption.

%

\vspace{2mm}

\ni\textbf{(HBQD)} There exists a constant $K_f$ such that the function $f$ satisfies
\beqs
| f(t,  x, y, z, u) - f(t',  x', y', z', u') | & \leq & K_{f} \big[ |x - x'| + |y - y'| + |u - u'|+ |t - t'|^{1\over 2} \big]  \\ & &+ L_{f, z} ( 1+|z| + |z'| ) |z - z'| \;,
\enqs
for all $(t,x,x',y,y',z,z',u,u')\in [0,T] \times \R^2\times \R^2\times \R^2\times \R^2$.\\

\vspace{2mm}

In the sequel of this section, the space of BMO martingales plays a key role for the a priori estimates of processes $Z^0$ and $Z^1$. We refer to \cite{kaz94} for the theory of BMO martingales. Here, we just give the definition of a BMO martingale and recall a property that we use in the sequel. 
\begin{Definition}
A process $M$ is said to be a $BMO_\F[0,T]$-martingale if 
 $M$ is a square integrable $\F$-martingale s.t. 
\beqs
{\|M\|}_{BMO_\F[0,T]} & := & \sup_{\tau \in \Tc_\F[0,T]} \E \Big[ \big|M_T-M_\tau \big|^2 \Big|\Fc_\tau \Big]^{1/2} ~ < ~ \infty \;,
\enqs
where $\Tc_\F[0,T]$ denotes the set of $\F$-stopping times valued in $[0,T]$.
\end{Definition}
\ni The BMO condition provides a property on the Dolean-Dade exponential of the process $M$. 
\begin{Lemma}\label{ppte-gen-BMO}
Let $M$  be a $BMO_\F[0,T]$-martingale. Then the stochastic exponential $\Ec(M)$ defined by 
\beqs
\Ec ( M)_t & = & \exp \Big( M_t - \frac{1}{2} \langle M,M \rangle\Big)_t \;, \quad 0 \leq t \leq T\;,
\enqs
is a uniformly integrable $\F$-martingale.
\end{Lemma}
\noindent We refer to \cite{kaz94} for the proof of this result.\\

We first state a BMO property for the processes $Z^0$ and $Z^1$, which will be used in the sequel to provide an estimate for these processes.
\begin{Lemma}\label{BMO}
Under \textbf{(HF)}, \textbf{(HBQ)} and \textbf{(HBQD)}, 
the martingales $\int_0^.Z_s^0 d W_s$ and  $\int_0^.Z^1_s(\theta) \mathds{1}_{s\geq\theta}d W_s$, $\theta\in[0,T]$ are $BMO_\F[0,T]$-martingales and there exists a constant $K$ which is independent from $\theta$ such that
\beqs
{\Big\|\int_0^.Z^0_s d W_s\Big\|}_{BMO_\F[0,T]} & \leq & K\;,\\
\sup_{\theta\in[0,T]}{\Big\|\int_0^.Z^1_s(\theta) \mathds{1}_{s\geq\theta}d W_s\Big\|}_{BMO_\F[0,T]} & \leq & K\;.
\enqs
 \end{Lemma}

\ni\textbf{Proof.} 
Define the function $\phi:~\R\rightarrow\R$ by
\beq\label{ppte phi}
\phi(x) &  = & (e^{2K_qx} - 2K_qx - 1)/2 K_q^2\;,\quad x\in\R\;.
\enq
 We notice that $\phi$ satisfies 
 \beqs
 \phi'(x)  & \geq &  0 ~\mbox{ and }~ \frac{1}{2} \phi''(x) - K_q \phi'(x) ~=~ 1\;,
 \enqs
 for $x \geq 0$. Since $Y^0$ and $Y^1(.)$ are solutions to quadratic BSDEs with bounded terminal conditions, we get from Proposition 2.1 in \cite{k00}  
 the existence of a constant $m$ such that 
 \beq\label{born-unif-Y}
\| Y^0\|_{\Sc^\infty[0,T]}  & \leq &  m \quad \mbox{ and }\quad  \sup_{\theta\in[0,T]}\|Y^1(\theta)\|_{\Sc^\infty[\theta,T]}  ~ \leq ~  m\;.
\enq
Applying It\^o's formula we get 
\begin{multline*}
\phi(Y^0_\nu + m) + \E \Big( \int_\nu^T \frac{1}{2} \phi''(Y^0_s + m) | Z^0_s|^2 ds\Big| \Fc_\nu \Big)  =  \\
\E(\phi( Y^0_T + m) | \Fc_\nu) + \E \Big( \int_\nu^T \phi'(Y^0_s + m) f(s, X^0_s, Y^0_s, Z^0_s, Y^1_s(s) - Y^0_s) ds \Big| \Fc_\nu \Big) \;.
\end{multline*}
for any $\F$-stopping time $\nu$ valued in $[0,T]$. 
From the growth assumption on the generator $f$ in \textbf{(HBD)}, \reff{ppte phi} and \reff{born-unif-Y}, we obtain
\begin{multline*}
\phi(Y^0_\nu + m) + \E \Big( \int_\nu^T| Z^0_s|^2 ds\Big| \Fc_\nu \Big)  \leq \\
 \E(\phi( Y^0_T + m) | \Fc_\nu)
 + \E \Big( \int_\nu^T \phi'(Y^0_s + m) K_q (1 + 2 {\|Y^0\|}_{\Sc^{\infty}} + \sup_{\theta\in[0,T]}\|Y^1 (\theta)\|_{\Sc^{\infty}[\theta,T]}) ds \Big| \Fc_\nu \Big) \;. 
\end{multline*}
This last inequality and \reff{born-unif-Y} imply that there exists a constant $K$ which depends only on $m$, $T$ and $K_q$ such that for all $\F$-stopping times $\nu \in[0, T]$
\beqs
 \E \Big( \int_\nu^T| Z^0_s|^2 ds\Big| \Fc_\nu \Big) & \leq & K \;.
\enqs
For the process $Z^1$, we use the same technics. Let us fix $\theta \in [0, T]$. 
Applying It\^o's fomula we get
\begin{multline*}
\phi(Y^1_{\nu\vee \theta}(\theta) + m) + \E \Big( \int_{\nu\vee \theta}^T \frac{1}{2} \phi''(Y^1_s(\theta) + m) | Z^1_s(\theta)|^2ds \Big| \Fc_{\nu \vee \theta} \Big)  =   \\
 \E(\phi( Y^1_T(\theta) + m) | \Fc_\nu) + \E \Big( \int_{\nu \vee \theta}^T \phi'(Y^1_s(\theta) + m) f(s, X^1_s(\theta), Y^1_s(\theta), Z^1_s(\theta), 0) ds \Big| \Fc_{\nu \vee \theta} \Big) \;,
\end{multline*}
for any $\F$-stopping time $\nu$ valued in $[0,T]$.
From the growth assumption on the generator $f$ in \textbf{(HBQ)}, \reff{ppte phi} and \reff{born-unif-Y}, we obtain
\begin{multline*}
\phi(Y^1_{\nu \vee \theta}(\theta) + m) + \E \Big( \int_{\nu \vee \theta}^T| Z^1_s(\theta)|^2ds \Big| \Fc_\nu \Big)   \leq   \E(\phi( Y^1_T(\theta) + m) | \Fc_\nu)\\
 + \E \Big( \int_{\nu \vee \theta}^T \phi'(Y^1_s(\theta) + m) K_q (1 + \|Y^1 (\theta)\|_{\Sc^{\infty}[\theta,T]}) ds \Big| \Fc_\nu \Big) \;.
\end{multline*}
This last inequality and \reff{born-unif-Y} imply that there exists a constant $K$ which depends only on $m$, $T$ and $K_q$, such that for all $\F$-stopping times $\nu$ valued in $[0,T]$
\beqs
 \E \Big( \int_\nu^T| Z^1_s(\theta)|^2\mathds{1}_{s\geq\theta}  ds\Big|\Fc_\nu \Big) & \leq & K \;.
\enqs
\ep

\subsection{Some bounds about $X^0$ and $X^1$}
In this part, we give some bounds about the processes $X^0$ and $X^1$ which are used to get a uniform bound for the processes $Z^0$ and $Z^1$.
\begin{Proposition}\label{borne sur X} Suppose that \textbf{(HF)} holds. Then, we have
 \beq 
| \nabla X^0_t |~ :=~ \Big| \frac{\partial X^0_t}{\partial x} \Big| & \leq & e^{L_a T} \;, \quad 0\leq t\leq T\;,\label{majX0}
\enq
and  for any $\theta \in [0,T]$ we have
 \beq 
| \nabla^\theta X^1_t(\theta) | ~:=~ \Big| \frac{\partial X^1_t(\theta)}{\partial X^1_\theta(\theta)} \Big|  & \leq & e^{L_a T} \;, \quad \theta\leq t\leq T\;,\label{majX1}
\enq
 \beq 
| \nabla X^1_t(\theta) |~ :=~\Big| \frac{\partial X^1_t(\theta)}{\partial x} \Big|& \leq &(1 + L_a e^{L_a T}) e^{L_a T} \;,  \quad \theta \leq t\leq T\;.\label{majX1-x}
\enq
\end{Proposition}
\ni\textbf{Proof.} 
We first suppose that $b$ and $\beta$ are $C^1_b$ w.r.t. $x$. By definition we have 
 \beqs
\nabla X^0_t & = & 1 + \int_0^t \partial_x b ( s, X^0_s ) \nabla X^0_s ds\;,  \quad 0 \leq t \leq T \;.
\enqs 
  We get from Gronwall's lemma
 \beqs 
| \nabla X^0_t | & \leq & e^{L_a T} \;, \quad \quad0\leq t\leq T\;.
\enqs
In the same way, we have
\beqs 
\nabla^\theta X^1_t(\theta) & = & 1 + \int_\theta^t \partial_x b ( s, X^1_s(\theta)) \nabla^\theta X^1_s(\theta) ds \;, \quad \quad \theta \leq t \leq T \;,
\enqs 
and from Gronwall's lemma we get
 \beqs 
| \nabla^\theta X^1_t | & \leq & e^{L_a T} \;, \quad \theta\leq t\leq T\;.
\enqs
Finally we prove the last inequality. By definition 
\beqs 
\nabla X^1_t(\theta)  & = & 1 + \int_0^t \partial_x b ( s, X^1_s(\theta)) \nabla X^1_s(\theta) ds + \partial_x \beta (\theta, X^0_\theta) \nabla X^0_\theta \;, \quad \quad \theta \leq t \leq T \;.
\enqs 
Using the inequality \reff{majX0}, we get
\beqs 
|\nabla X^1_t(\theta)| & \leq & 1+ L_a e^{L_a T} + \int_0^t L_a |\nabla X^1_s(\theta)| ds  \;, \quad \quad \theta \leq t \leq T \;,
\enqs 
from Gronwall's lemma we get
 \beqs 
| \nabla X^1_t(\theta) | & \leq &(1 + L_a e^{L_a T}) e^{L_a T} \;,  \quad \theta \leq t\leq T\;.
\enqs
When $b$ and $\beta$ are not differentiable, we can also prove the result by regularization. We consider a density $q$ which is $C^\infty_b$ on $\R$ with a compact support, and we define an approximation $(b^\epsilon,\beta^\epsilon)$ of $(b, \beta)$ in $C^1_b$ by
\beqs
(b^\epsilon,  \beta^\epsilon)(t, x) & = & \frac{1}{\eps}\int_{\R} (b, \beta) (t, x') q \Big( \frac{x-x'}{\epsilon} \Big) dx'\;,\quad (t,x)\in[0,T]\times\R\;.
\enqs
We then use the convergence of $(X^{0,\epsilon},X^{1, \epsilon} (\theta))$ to $(X^0,X^{1} (\theta))$ and we get the result.
\ep

\subsection{Some bounds about $Y^0$ and $Y^1$}
In this part, we give some bounds about the processes $Y^0$ and $Y^1$ which are used to get a uniform bound for the processes $Z^0$ and $Z^1$.
\begin{Lemma}
Suppose that \textbf{(HF)}, \textbf{(HBQ)} and \textbf{(HBQD)} hold. Then, for any $\theta \in [0, T]$
\beq \label{majY1}
| \nabla^\theta Y^1_t(\theta)| ~:= ~\Big| \frac{\partial Y^1_t(\theta)}{\partial X^1_\theta(\theta)} \Big| & \leq & e^{( L_a + K_f ) T } ( K_g + T K_f ) \;,\quad \theta\leq t\leq T\;.
\enq 
\end{Lemma}

\ni\textbf{Proof.}
We first suppose that $b$, $f$ and $g$ are $C^1_b$ w.r.t. $x$, $y$ and $z$. In this case 
$(X^1(\theta) ,  Y^1(\theta), Z^1(\theta))$ is also differentiable w.r.t. $X^1_\theta(\theta) $ and we have
\beq \label{gradient de Y}
 \nabla^\theta Y^1_t(\theta) 
 & = & \nabla g ( X^1_T(\theta)  ) \nabla^\theta X^1_T(\theta)  - \int_t^T \nabla^\theta Z^1_s(\theta) dW_s \\
& & \nonumber + \int_t^T  \nabla f (s, X^1_s(\theta) , Y^1_s(\theta), Z^1_s(\theta) , 0 )  \big( \nabla^\theta X^1_s(\theta) ,\nabla^\theta Y^1_s(\theta) , \nabla^\theta Z^1_s(\theta) \big)  ds \;,
\enq
for $t\in[\theta,T]$. 
Define the process $R(\theta)$ by
\beqs
R_t(\theta) & := & \exp\Big({\int_0^t \partial_y  f(s, X^1_s(\theta), Y^1_s(\theta), Z^1_s(\theta),0)\mathds{1}_{s\geq \theta}ds}\Big)\;,\quad 0\leq t\leq T\;.
\enqs  
\ni Applying It\^o's formula, we get
\beq\nonumber
R_t(\theta) \nabla^\theta Y^1_t(\theta) & = & R_T(\theta) \nabla g(X^1_T(\theta)) \nabla^\theta X^1_T(\theta)  \\
& &  + \int_t^T R_s(\theta) \partial_x f (s, X^1_s(\theta), Y^1_s(\theta), Z^1_s(\theta),0 ) \nabla^\theta X^1_s(\theta) d s \nonumber\\
& & - \int_t^T R_s(\theta) \nabla^\theta Z^1_s(\theta) d W^1_s(\theta) \;, \quad \theta \leq t \leq T \;,\label{itoY1}
\enq 
where  the process $W^1(\theta)$ is defined by
\beq\label{W1}
 W^1_t(\theta) := W_t- \int_0^t\partial_z f (s, X^1_s(\theta), Y^1_s(\theta), Z^1_s(\theta), 0)\mathds{1}_{s\geq\theta} ds
\enq
for $t\in[0,T]$.
 From \textbf{(HBQD)}, there exists a constant $K >0$ such that we have 
\beqs
\Big\| \int_0^.\partial_z f (s, X^1_s(\theta), Y^1_s(\theta), Z^1_s(\theta), 0)\mathds{1}_{s\geq\theta} dW_s \Big\|^2_{BMO_\F[0,T]} & \leq  &\\
 K  \Big( 1+\sup_{\vartheta\in\Tc_\F[0,T]} \E\Big[\int_\vartheta^T|Z_s^1(\theta)|^2\mathds{1}_{s \geq \theta} ds\Big|\Fc_\vartheta\Big]\Big) 
 & \leq & \\
 K \Big(1+\Big\|\int_0^.Z^1_s(\theta)\mathds{1}_{ s\geq \theta}dW_s\Big\|^2_{BMO_\F[0,T]}\Big)
  & < & \infty\;,
\enqs
where the last inequality comes from Lemma \ref{BMO}. 

Hence by Lemma \ref{ppte-gen-BMO} the process $\Ec(\int_0^.\partial_z f (s, X^1_s(\theta), Y^1_s(\theta), Z^1_s(\theta), 0) \mathds{1}_{s\geq\theta}dW_s )$ is a uniformly integrable martingale.
Therefore, under the probability measure $\Q^1(\theta)$ defined by 
\beqs
\frac{d \Q^1(\theta)}{d \P} \Big|_{\Fc_t} & := & \Ec\Big(\int_0^.\partial_z f (s, X^1_s(\theta), Y^1_s(\theta), Z^1_s(\theta), 0) \mathds{1}_{s\geq\theta}dW_s \Big)_t\;,\quad 0\leq t\leq T\;,
\enqs
 we can apply Girsanov's theorem and $W ^1(\theta)$ is a Brownian motion under the probability measure $\Q^1(\theta)$. We then get from  \reff{itoY1}
\beqs
R_t(\theta) \nabla^\theta Y^1_t & = & \E_{\Q^1(\theta)}\Big[R_T(\theta) \nabla g(X^1_T) \nabla^\theta X^1_T  
 + \int_t^T R_s(\theta) \partial_x f (s, X^1_s, Y^1_s, Z^1_s,0 ) \nabla^\theta X^1_s d s \Big|\Fc_t\Big]\;.
\enqs
This last equality, \textbf{(HBQD)} and \reff{majX1}  give
\beq \label{majY1}
| \nabla^\theta Y^1_t(\theta)| & \leq & e^{( L_a + K_f ) T } ( K_g + T K_f ) \;,\quad \theta\leq t\leq T\;.
\enq 
 When $b$, $f$ and $g$ are not $C^1_b$, we can also prove the result by regularization as for Proposition \ref{borne sur X}.
\ep

\begin{Lemma}Suppose that \textbf{(HF)},  \textbf{(HBQ)} and \textbf{(HBQD)} hold. Then,
\beq
| \nabla Y^1_t(t)|  &\leq&  (1+L_a e^{L_a T})e^{( L_a + K_f ) T } ( K_g + T K_f )\;,  \quad 0 \leq t \leq T \;. \label{maj nabla Y1}
\enq
\end{Lemma}

 \ni\textbf{Proof.}
 Firstly, we suppose that $b$, $\beta$, $g$ and $f$ are $C^1_b$ w.r.t. $x$, $y$ and $z$. Then, for any $t \in [0,T]$
 \beq \label{gradient de Y1}
 \nabla Y^1_t(t) ~:= ~ \frac{\partial Y^1_t(t)}{\partial x } & = & \nabla g ( X^1_T(t)) \nabla X^1_T(t)\\
\nonumber & & + \int_t^T  \nabla f (s, X^1_s(t), Y^1_s(t), Z^1_s(t), 0)(\nabla X^1_s(t),\nabla Y^1_s(t),\nabla Z^1_s(t))ds\\
 & &  \nonumber   - \int_t^T \nabla Z^1_s(t) dW_s\;.
\enq
\ni Applying It\^o's formula, we get
\beq\nonumber
R_t(t) \nabla Y^1_t(t) & = & R_T(t) \nabla g(X^1_T(t)) \nabla X^1_T(t)  \\
& &  + \int_t^T R_s(t) \partial_x f (s, X^1_s(t), Y^1_s(t), Z^1_s(t),0 ) \nabla X^1_s(t) d s \nonumber\\
& & - \int_t^T R_s(t) \nabla Z^1_s(t) d W^1_s(t) \;, \quad 0 \leq t \leq T \;,\label{itoY1-x}
\enq 
where  the process $W^1(.)$ is defined in \reff{W1}. We have proved previously that  $W ^1(t)$ is a Brownian motion under the probability measure $\Q^1(t)$. We then get 
\beqs
R_t(t) \nabla Y^1_t(t) & = & \E_{\Q^1(t)}\Big[R_T(t) \nabla g(X^1_T(t)) \nabla X^1_T(t)  
 + \int_t^T R_s(t) \partial_x f (s, X^1_s(t), Y^1_s(t), Z^1_s(t),0 ) \nabla X^1_s(t) d s \Big|\Fc_t\Big]\;.
\enqs
This last inequality, \textbf{(HBQD)} and \reff{majX1}  give
\beqs
| \nabla Y^1_t(t)|  &\leq&  (1+L_a e^{L_a T})e^{( L_a + K_f ) T } ( K_g + T K_f )\;,  \quad 0 \leq t \leq T \;.
\enqs
When $b$, $f$ and $g$ are not $C^1_b$, we can also prove the result by regularization as for Proposition \ref{borne sur X}.
 \ep
 
\begin{Lemma}\label{nabla Y 0}
Suppose that \textbf{(HF)},  \textbf{(HBQ)} and \textbf{(HBQD)}  hold. Then, 
\beqs
| \nabla Y^0_t | & \leq & e^{ (2 K_f + L_a) T } (K_g + K_f T) \big( 1 + T K_f e^{K_fT} (1 + L_a e^{L_a T})\big) \;, \quad 0 \leq t \leq T\;.
\enqs 
\end{Lemma}
 
 \ni\textbf{Proof.}
We first suppose that $b$, $\beta$, $g$ and $f$ are $C^1_b$ w.r.t. $x$, $y$, $z$ and $u$, then $(X^0, Y^0, Z^0)$ is differentiable w.r.t. $x$ and  we have 
 \beq \label{gradient de Y0}
\nabla Y^0_t& = & \nabla g ( X^0_T) \nabla X^0_T\\
\nonumber & & + \int_t^T \Big( \nabla f (s, X^0_s, Y^0_s, Z^0_s, Y^1_s(s) - Y^0_s)(\nabla X^0_s,\nabla Y^0_s,\nabla Z^0_s,\nabla Y^1_s(s)-\nabla Y^0_s)ds\\
  & &  \nonumber 
  - \int_t^T \nabla Z^0_s dW_s\;,  \quad 0 \leq t \leq T\;.
\enq
Define the process $R^0$ by 
\beqs
R_t^0 & := & \exp\Big(\int_0^t (\partial_y - \partial_u) f (s, X^0_s, Y^0_s, Z^0_s, Y^1_s(s) - Y^0_s)ds\Big)\;,\quad 0\leq t\leq T\;.
\enqs
Applying It\^o's fomula we have
\beqs
R^0_t\nabla Y^0_t & = & R_T^0 \nabla g(X^0_T) \nabla X^0_T  \\
& &  + \int_t^T R_s^0 \partial_x f (s, X^0_s, Y^0_s, Z^0_s, Y^1_s(s) - Y^0_s ) \nabla X^0_s d s \\
& & + \int_t^T R_s^0   \partial_u f (s, X^0_s, Y^0_s, Z^0_s, Y^1_s(s) - Y^0_s ) \nabla Y^1_s ( s )  ds\\
& & - \int_t^T R_s^0  \nabla Z^0_s d {W}^0_s 
\enqs
where $d {W}^0_s := d W_s - \partial_z f (s, X^0_s, Y^0_s, Z^0_s, Y^1_s(s) - Y^0_s) ds$. From \textbf{(HBQD)}, there exists a constant $K >0$ such that we have 
\beqs
\Big\| \int_0^.\partial_z f (s, X^0_s, Y^0_s, Z^0_s, Y^1_s(s) - Y^0_s) dW_s \Big\|^2_{BMO_\F[0,T]} & \leq  & \\K  \Big( 1+\sup_{\vartheta\in\Tc_\F[0,T]} \E\Big[\int_\vartheta^T|Z_s^0|^2 ds\Big|\Fc_\vartheta\Big]\Big)
  &\leq& \\  K \Big(1+\Big\|\int_0^.Z^0_sdW_s\Big\|^2_{BMO_\F[0,T]}\Big)
  & < & \infty\;,
\enqs
where the last inequality comes from Lemma \ref{BMO}. 

Hence by Lemma \ref{ppte-gen-BMO}  the process $\Ec(\int_0^.\partial_z f (s, X^0_s, Y^0_s, Z^0_s, Y^1_s(s) - Y^0_s) dW_s )$ is a uniformly integrable martingale. Therefore, under the probability measure $\Q^0$ defined by
\beqs
\frac{d \Q^0}{d \P} \Big|_{\Fc_t} &  := &  \Ec\Big(\int_0^.\partial_z f (s, X^0_s, Y^0_s, Z^0_s, Y^1_s(s) - Y^0_s) dW_s\Big)_t
\enqs
 we can apply Girsanov's theorem and $W ^0$ is a Brownian motion under the probability measure $\Q^0$. Then, we get 
\beqs
R_t^0 \nabla Y^0_t & = & \E_{\Q^0} \Big[ R_T^0 \nabla g(X^0_T) \nabla X^0_T  \\
& &  + \int_t^T R_s^0 \partial_x f (s, X^0_s, Y^0_s, Z^0_s, Y^1_s(s) - Y^0_s ) \nabla X^0_s d s \\
& & + \int_t^T R_s^0 \partial_u f (s, X^0_s, Y^0_s, Z^0_s , Y^1_s(s) - Y^0_s) \nabla Y^1_s ( s )  ds \Big| \Fc_t \Big] \;.
\enqs
Using inequalities \reff{majX0} and \reff{maj nabla Y1} we get
\beqs
| \nabla Y^0_t | & \leq & e^{ (2 K_f + L_a) T } (K_g + K_f T) \big( 1 +  K_f T e^{K_fT} (1 + L_a e^{L_a T})\big) \;, \quad 0 \leq t \leq T\;.
\enqs
When $b$, $\beta$, $f$ and $g$ are not $C^1_b$, we can also prove the result by regularization as for Proposition \ref{borne sur X}.
\ep

\subsection{A uniform bound for $Z^0$ and $Z^1$}
Using the previous bounds, we obtain a uniform bound for the processes $Z^0$ and $Z^1$.
\begin{Proposition} \label{borne 1 pour Z 1}
Suppose that \textbf{(HF)}, \textbf{(HBQ)} and \textbf{(HBQD)} 
  hold. Then, for any $\theta \in [0, T]$, there exists a version of $Z^1 ( \theta )$ such that
\beqs
| Z^1_t ( \theta ) | & \leq &e^{(  2L_a + K_f ) T } ( K_g + T K_f ) K_a \;, \quad \theta \leq t \leq T \;.
\enqs 
\end{Proposition}

\ni\textbf{Proof.} 
Using Malliavin calculus, we have the classical representation of the process $Z^1(\theta)$ given by $\nabla^\theta Y^1(\theta)(\nabla^\theta X^1(\theta))^{-1}\sigma(.)$ (see \cite{kl11a}). In the case where $b$, $f$ and $g$ are $C^1_b$ w.r.t. $x$, $y$ and $z$, we obtain from \reff{majY1}
\beqs 
| Z^1_t(\theta)  | & \leq &  e^{(  2L_a + K_f ) T } ( K_g + T K_f ) K_a \quad a.s.
\enqs
since $|(\nabla^\theta X^1(\theta))^{-1}| \leq e^{L_a T}$ (the proof of this inequality is similar to the one of  \reff{majX1}).\\
When $b$, $f$ and $g$ are not differentiable, we can also prove the result by a standard approximation and stability results for BSDEs with linear growth. 
\ep

\vspace{3mm}

\begin{Proposition}\label{borne 1 pour Z 0}
Suppose that \textbf{(HF)},  \textbf{(HBQ)} and \textbf{(HBQD)} hold. Then, there exists a version of $Z^0$ such that
\beqs
| Z^0_t | &\leq &e^{ 2 ( K_f + L_a) T } (K_g + K_f T) \big( 1 + T K_f e^{K_fT} (1 + L_a e^{L_a T})\big) K_a \;,  \quad 0 \leq t \leq T \;.
\enqs 
\end{Proposition}
 
 \ni\textbf{Proof.}
Thanks to the Malliavin calculus, it is classical to show that a version of $Z^0$ is given by $\nabla Y^0 (\nabla X^0)^{-1} \sigma(.)$ (see \cite{kl11a}). So, in the case where $b$, $\beta$, $g$ and $f$ are $C^1_b$ w.r.t. $x$, $y$, $z$ and $u$, we obtain from \reff{majX0} and Lemma \ref{nabla Y 0}
\beqs
| Z^0_t | & \leq & e^{ 2 ( K_f + L_a) T } (K_g + K_f T) \big( 1 + T K_f e^{K_fT} (1 + L_a e^{L_a T})\big) M_\sigma \quad a.s. 
\enqs
since $| (\nabla X^0_t)^{-1} | \leq e^{L_a T}$ (the proof of this inequality is similar to the one of \reff{majX0}).

When $b$, $\beta$, $g$ and $f$ are not differentiable, we can also prove the result by a standard approximation and stability results for BSDEs with linear growth. 
 \ep
 
\section{Discrete-time approximation for the FBSDE}

 \subsection{Discrete-time scheme for the FBSDE}
 Throughout the sequel, we consider a discretization grid $\pi :=\{t_0,\ldots,t_n\}$ of $[0,T]$ with $0=t_0<t_1<\ldots<t_n=T$.  For $t\in[0,T]$, we denote by $\pi(t)$ the largest element of $\pi$ smaller than $t$
\beqs
\pi(t) & := & \max\big\{\;t_i\;,~i=0,\ldots,n~|~t_i\leq t\;\big\}\;.
\enqs
We also denote by $|\pi|$ the mesh of $\pi$
\beqs
|\pi| & := & \max\big\{\;t_{i+1}-t_i\;,~i=0,\ldots,n-1\;\big\}\;,
\enqs
that we suppose satisfying $|\pi|\leq 1$,  and by $\Delta W^\pi_i$ (resp. $\Delta t^\pi_i $) the increment of $W$ (resp. the difference) between $t_i$ and $t_{i-1}$: $\Delta W^\pi_i := W_{t_{i}} - W_{t_{i-1}}$ (resp. $\Delta t^\pi_i :=t_i-t_{i-1}$), 
for $1\leq i\leq n$.

\vspace{2mm}

We introduce an approximation of the process $X$ based on the discretization of the processes $X^0$ and $X^1$.

\vspace{2mm}

\ni$\bullet$ \emph{Euler scheme for $X^0$}. We consider the classical scheme $X^{0,\pi}$ defined by
\begin{equation}\label{Euler0} \left\{
\begin{aligned}
  \XzePi_{t_0} ~ & = ~ x \;,\\
  \XzePi_{t_i} ~ & = ~   \XzePi_{t_{i-1}} + b(t_{i-1}, \XzePi_{t_{i-1}})\Delta t_i^\pi + \sigma(t_{i-1}) \Delta W^\pi_i  \;,\quad 1\leq i\leq n\;.\\
   \end{aligned}
\right.\end{equation}
\ni$\bullet$ \emph{Euler scheme for $X^1$}. Since the process $X^1$ depends on two parameters $t$ and $\theta$, we introduce a discretization of $X^1$ in these two variables. We then consider the following scheme\footnote{$\delta_{i=k} = 0$ if $i \neq k$ and $\delta_{i=k} = 1$ if $i=k$.}
\begin{equation}\label{Euler1} \left\{
\begin{aligned}
  \XunPi_{t_0}(t_k) ~ & = ~ x+\beta(t_0,x)\delta_{k=0} \;,\quad 0\leq k\leq N\;,\\
  \XunPi_{t_i}(t_k) ~ & = ~   \XunPi_{t_{i-1}}(t_k) + b(t_{i-1}, \XunPi_{t_{i-1}}(t_k))\Delta t^\pi_i + \sigma( t_{i-1}) \Delta W^\pi_i\\
   & \quad  + \beta(t_{i-1}, \XunPi_{t_{i-1}}(t_k)) \delta_{i=k} \;,\quad\ 1\leq i,k\leq n\;.\\
   \end{aligned}
\right.\end{equation}

\ni We are now able to provide an approximation of the process $X$ solution to the FSDE \reff{FSDE}.
We consider the scheme $X^\pi$ defined by 
\beq\label{EulerX}
X^\pi_{t} & = & \XzePi_{\pi(t)}\1_{t < \tau} +\XunPi_{\pi(t)}(\pi(\tau))\1_{t \geq \tau}\;,\quad \quad 0 \leq t \leq T\;.
\enq
We shall denote by $\{\Fc^{0,\pi}_i\}_{0 \leq i \leq n}$ (resp. $\{\Fc^{1,\pi}_i(\theta)\}_{0 \leq i \leq n}$) the discrete-time filtration  associated with $X^{0,\pi}$ (resp. $X^{1,\pi}(\theta)$)
\beqs
\Fc^{0,\pi}_i& := & \sigma( \XzePi_{t_j},\; j \leq i) \\
 \mbox{ (resp. }\Fc^{1,\pi}_i(\theta)  & := &  \sigma( \XunPi_{t_j}(\theta),\; j \leq i)\mbox{)}\;.%
\enqs

\vspace{2mm}

We introduce an approximation of $(Y,Z)$ based on the discretization of $(Y^0,Z^0)$ and $(Y^1,Z^1)$.   
To this end we introduce the backward implicit schemes on $\pi$ associated with the BSDEs \reff{BSDE1} and \reff{BSDE0}. Since the system is recursively coupled, we first introduce the scheme associated with  \reff{BSDE1}. We then use it to define the scheme associated with \reff{BSDE0}. 

\vspace{2mm}

\ni$\bullet$   \textit{Backward Euler scheme for $(Y^1,Z^1)$}. We consider the classical implicit scheme $(\YunPi,\ZunPi)$ defined by 
\begin{equation}\label{YZ1} \left\{
\begin{aligned}
  \YunPi_{T}(\pi ( \theta )) ~  = & ~ g ( \XunPi_{T}(\pi ( \theta )) ) \;,\\
  \YunPi_{t_{i-1}}(\pi ( \theta )) ~ = & ~ \E \big[ \YunPi_{t_{i}}(\pi ( \theta )) \big| \Fc^{1,\pi}_{i-1}(\pi(\theta)) \big] \\
   & + f ( t_{i-1},\XunPi_{t_{i-1}}(\pi ( \theta )), \YunPi_{t_{i-1}}(\pi ( \theta )), \ZunPi_{t_{i-1}}(\pi ( \theta )), 0) \Delta t_i^\pi \;,\\
    \ZunPi_{t_{i-1}}(\pi ( \theta )) ~ = & ~ \frac{1}{\Delta t_i^\pi } \E \big[ \YunPi_{t_i}(\pi ( \theta )) \Delta W^\pi_{i}  \big| \Fc^{1,\pi}_{i-1}(\pi(\theta)) \big] \;,\quad t_{i-1}\geq \pi(\theta)\;.
    \end{aligned}
\right.\end{equation}
  \vspace{2mm}

\ni$\bullet$  \textit{Backward Euler scheme for $(Y^0,Z^0)$}. Since the generator of \reff{BSDE0} involves the process ${(Y^1_t(t))}_{t\in[0,T]}$, we consider a discretization based on $\YunPi$. We therefore consider the scheme $(\YzePi,\ZzePi)$ defined by 
\begin{equation}\label{YZ0} \left\{
\begin{aligned}
  \YzePi_{T} ~ & = ~ g ( \XzePi_{T} ) \;,\\
  \YzePi_{t_{i-1}} ~ & = ~ \E\big[ \YzePi_{t_{i}} \big| \Fc^{0,\pi}_{i-1} \big] + \bar f^\pi ( t_{i-1},\XzePi_{t_{i-1}}, \YzePi_{t_{i-1}}, \ZzePi_{t_{i-1}}) \Delta t^\pi_{i} \;,\\
    \ZzePi_{t_{i-1}} ~ & = ~ \frac{1}{\Delta t^\pi_{i} } \E \big[ \YzePi_{t_i} \Delta W^\pi_{i} \big| \Fc^{0,\pi}_{i-1}\big] \;,\quad1\leq i\leq n\;,
    \end{aligned}
\right.\end{equation}
  where $\bar f^\pi(t,x,y,z) : =  f\big(t,x,y,z,\YunPi_{\pi(t)}(\pi(t))-y\big)$
  for all $(t,x,y,z)\in[0,T]\times\R \times\R\times\R$.
  
\ni   We then consider the following scheme defined by
\begin{equation}\label{SchemeYZ} \left\{
\begin{aligned}
Y^\pi_{t} ~& =~  Y^{0, \pi}_{\pi(t)} \mathds{1}_{t < \tau} + \YunPi_{\pi(t)} ( \pi ( \tau ) ) \mathds{1}_{t \geq  \tau} \,,\\
Z^\pi_{t}~ & =~  Z^{0, \pi}_{\pi(t)} \mathds{1}_{t \leq  \tau} + \ZunPi_{\pi(t)} ( \pi (  \tau ) ) \mathds{1}_{t >  \tau} \,,\\
U^\pi_{t} ~& =~  \big( \YunPi_{\pi(t)} ( \pi(t)) - Y^{0, \pi}_{\pi(t)} \big) \mathds{1}_{t \leq  \tau} \,,
    \end{aligned}
\right.\end{equation}
for $t\in[0,T]$.

\subsection{Convergence of the scheme for the FBSDE}

We now concentrate on the error approximation of the processes $X$, $Y$, $Z$ and $U$ by their scheme $X^\pi$, $Y^\pi$, $Z^\pi$ and $U^\pi$. To this end we introduce extra assumptions on the regularity of the forward coefficients w.r.t. the time variable $t$.

\vspace{2mm}

\ni\textbf{(HFD)} There exists a constant $K_t$ such that the functions $b$, $\sigma$ and $\beta$ satisfy
\beqs
\big|b(t,x)-b(t',x)\big| +\big|\sigma(t)-\sigma(t')\big| & \leq & K_t|t-t'|^{1\over 2}\;,\\
\big|\beta(t,x)-\beta(t',x)\big| +\big|\sigma(t)-\sigma(t')\big| & \leq & K_t|t-t'|\;,
\enqs
for all $(t,t',x)\in[0,T]^2 \times\R$.

\vspace{2mm}

\ni We can now state the main result of the paper.

\begin{Theorem}
Under  \textbf{(HF)}, \textbf{(HFD)}, \textbf{(HBQ)} and \textbf{(HBQD)}  we have the following estimate
\beqs
\E \Big[ \sup_{t\in[0,T]} \big| X_{t} - X^\pi_{t} \big|^2 \Big] + \sup_{t\in[0,T] } \E \Big[ \big| Y_{t} - Y^{\pi}_{t} \big|^2 \Big]  & & \\
+ \E \Big[ \int_0^T \big| Z_t - Z^\pi_{t} \big|^2 dt \Big] + \E \Big[ \int_0^T \lambda_t \big| U_t - U^\pi_{t} \big|^2   dt \Big] & \leq & K |\pi| \;,
\enqs
for a constant $K$ which does not depend on $\pi$.
\end{Theorem}
\ni\textbf{Proof.} Fix $M\in\R$ such that  
\beqs
M & \geq &  \max\Big\{e^{(  2L_a + K_f ) T } ( K_g + T K_f ) K_a ~; \\
 & &   \quad \quad e^{ 2 ( K_f + L_a) T } (K_g + K_f T) \big( 1 + T K_f e^{K_fT} (1 + L_a e^{L_a T})\big) K_a \Big\}\;,
\enqs 
and define the function $\tilde f$ by 
\beqs
\tilde f (t,x,y,z,u) & = & f (t,x,y,\varphi_M(z),u) \;,\quad (t,x,y,z,u)\in[0,T]\times\R\times\R\times\R\times\R\;,
\enqs 
 where 
 \begin{equation*}
 \varphi_M(z)    := 
 \left\{
 \begin{array}{ccc}
 z & \mbox{ if }& |z|\leq M \\
M\frac{z}{|z|}    & \mbox{ if }& |z|> M
 \end{array}
 \right. , \quad z\in\R\;.
 \end{equation*}
 We notice that $\varphi_M$ is Lipschitz continuous and bounded. Therefore  we obtain from  $\textbf{(HBQD)}$ that $\tilde f$ is Lipschitz continuous.  
 
 Moreover, using Propositions \ref{borne 1 pour Z 1} and \ref{borne 1 pour Z 0}, we get that under 
 \textbf{(HF)}, \textbf{(HBQ)} and \textbf{(HBQD)}, 
 $(X,Y,Z)$ is also solution to the Lipschitz FBSDE
 \beqs
  X_t  & = & x + \int_0^t b(s, X_s) ds + \int_0^t \sigma(s) dW_s + \int_0^t \beta(s, X_{s^-} ) dH_s \;, \quad 0 \leq t\leq T\;,\\ \nonumber
  Y_t  & =  & g(X_T) + \int_t^T\tilde  f\big(s, X_s, Y_s, Z_s, U_s(1-H_s)\big) ds\\ 
   & & \quad \quad \quad \quad  - \int_t^T Z_s dW_s - \int_t^T U_s dH_s\;, \quad 0 \leq t\leq T\;. ~~
 \enqs
 Applying Theorem 5.1 in \cite{kl11a}, we get the result.
\ep

\end{document}